\documentclass{elsart}

\usepackage{latexsym}
\usepackage[psamsfonts]{amssymb}
\usepackage{amsfonts}
\usepackage{graphicx,amssymb,cite}
\usepackage{amsmath}
\usepackage{amssymb}
\usepackage[mathscr]{eucal}
\usepackage{enumerate}

\begin{document}

\newtheorem{tm}{Theorem}[section]
\newtheorem{pp}{Proposition}[section]
\newtheorem{lm}{Lemma}[section]
\newtheorem{df}{Definition}[section]
\newtheorem{tl}{Corollary}[section]
\newtheorem{re}{Remark}[section]
\newtheorem{eap}{Example}[section]

\newcommand{\pof}{\noindent {\bf Proof} }
\newcommand{\ep}{$\quad \Box$}

\newcommand{\al}{\alpha}
\newcommand{\be}{\beta}
\newcommand{\var}{\varepsilon}
\newcommand{\la}{\lambda}
\newcommand{\de}{\delta}
\newcommand{\st}{\stackrel}

\allowdisplaybreaks

\begin{frontmatter}

\title{Characterizations of precompact sets in fuzzy star-shaped numbers space with $L_p$-metric}
 \author{Huan Huang}
\author{}\ead{   hhuangjy@126.com   }
\address{ Department of Mathematics, Jimei
University, Xiamen 361021, China}

\date{}
\maketitle

\end{frontmatter}

\section{Introduction}

Compactness of fuzzy sets attracts many authors' attention.
The readers may see Ref. \cite{da,fan,greco,huang,ma,roman,wu2} for some details.
In
 this paper, we consider the characterizations of precompact sets in fuzzy star-shaped number space with $L_p$-metric.

\section{Preliminaries}

In this section,
we introduce some basic notations and results.
 For details, we refer the readers to references
\cite{wu, da}.

Let $\mathbb{N}$ be the set of all natural numbers, $\mathbb{R}^m$
be $m$-dimensional Euclidean space,
$K(  \mathbb{R}^m  )$
be
the
set of all nonempty compact set in $\mathbb{R}^m$,
and
$C(  \mathbb{R}^m  )$
be
the
set of all nonempty closed set in $\mathbb{R}^m$.
The
 well-known
 {\rm Hausdorff} metric $H$ on
   $C(\mathbb{R}^m)$ is defined by:
$$H(U,V)=\max\{H^{*}(U,V),\ H^{*}(V,U)\}$$
for arbitrary $U,V\in K(\mathbb{R}^m)$, where
  $$H^{*}(U,V)=\sup\limits_{u\in U}\,d\, (u,V) =\sup\limits_{u\in U}\inf\limits_{v\in
V}d\, (u,v).$$

A set $K\in K(\mathbb{R}^m)$ is said to be star-shaped relative to a point $x\in K$
if
for each $y\in K$, the line $\overline{xy}$ joining $x$ to $y$
is contained in $K$.
The kernel ker $K$ of $K$ is the set of all points $x\in K$
such that
$\overline{xy}\subset K$
for each $y\in K$.
The symbol $K_S(\mathbb{R}^m)$
is used to
denote
all the star-shaped sets in $\mathbb{R}^m$.

$K_C(\mathbb{R}^m)$
is used to
denote
all the nonempty compact and convex sets in $\mathbb{R}^m$.
Obviously,
$K_C(\mathbb{R}^m)\varsubsetneqq K_S(\mathbb{R}^m)$.

We use $F(\mathbb{R}^m)$ to represent all
fuzzy subsets on $\mathbb{R}^m$, i.e. functions from $\mathbb{R}^m$
to $[0,1]$.
For
$u\in F(\mathbb{R}^m)$, let $[u]_{\al}$ denote the $\al$-cut of
$u$, i.e.
\[
[u]_{\al}=\begin{cases} \{x\in \mathbb{R}^m : u(x)\geq \al \},\
\al\in(0,1],\\ {\rm supp}\, u=\overline{\{x \in
\mathbb{R}^m: u(x)>0\}}, \al=0.
\end{cases}
\]
For
$u\in F(\mathbb{R}^m)$,
we suppose that
\\
(\romannumeral1) \ $u$ is normal: there exists at least one $x_{0}\in \mathbb{R}^m$
with $u(x_{0})=1$;
\\
(\romannumeral2) \ $u$ is upper semi-continuous;
\\
(\romannumeral3) \ $u$ is fuzzy convex: $u(\la x+(1-\la)y)\geq {\rm min} \{u(x),u(y)\}$
for $x,y \in \mathbb{R}^m$ and $\la \in [0,1];$
\\
(\romannumeral4) \ $[u]_\lambda$ is a star-shaped set for all $\lambda\in (0,1]$.
\\
(\romannumeral5) \ $[u]_0$ is a compact set in $\mathbb{R}^m$.
\\
(\romannumeral6) \ $\left (\int_0^1   H( [u]_\al, \{0\}   )^p  \; d\alpha   \right)^{1/p}  <   +\infty $,
where $p\geq 1$ and $0$ denotes the origin of $\mathbb{R}^m$.

\begin{itemize}
  \item
  If
     $u$ satisfies (\romannumeral1), (\romannumeral2), (\romannumeral3) and (\romannumeral5),
  then
  we say $u$ is a fuzzy number. The set of all fuzzy numbers is denoted by $E^m$.

  \item
  If $u$ satisfies (\romannumeral1), (\romannumeral2), (\romannumeral4) and (\romannumeral5),
  then
  we say $u$ is a fuzzy  star-shaped number. The set of all fuzzy star-shaped  numbers is denoted by
  $S^m$.
\end{itemize}

Clearly, $E^m\subsetneq S^m$. Diamond \cite{da} introduced
the
 $d_p$ distance ($1\leq
p<\infty$) on $S^m$
which are defined by
\begin{equation}\label{dpm}
d_p\, (u,v)=\left(     \int_0^1 H([u]_\al, [v]_\al) ^p  \; d\alpha  \right)^{1/p},
\end{equation}
for all $u,v\in S^m$.
They pointed out that $d_p$ is a metric on $S^m$.

Wu and Zhao \cite{wu2} considered the compactness criteria of $E^m$ and $S^m$ with $d_p$ metric.

If we relax assumption (\romannumeral5) a little to obtain assumption (\romannumeral6),
then we obtain
 the following
  $L_p$--type noncompact fuzzy numbers or fuzzy star-shaped numbers.

\begin{df}
{\rm
If
  $u\in  F(\mathbb{R}^m)$ satisfies (\romannumeral1), (\romannumeral2), (\romannumeral3) and (\romannumeral6),
then
  we say $u$ is a $L_p$-type non-compact fuzzy number. The set of all such fuzzy numbers is denoted by $\widetilde{E}^{m,p}$.
  If there is no confusion, we also denote it by $\widetilde{E}^m$ for simplicity.

  If
  $u\in  F(\mathbb{R}^m)$ satisfies (\romannumeral1), (\romannumeral2), (\romannumeral4) and (\romannumeral6)
then
  we say $u$ is a $L_p$-type non-compact fuzzy  star-shaped number. The set of all such fuzzy  star-shaped numbers is denoted by $\widetilde{S}^{m,p}$.
  If there is no confusion, we also denote it by $\widetilde{S}^m$ for writing convenience.
  }
\end{df}

Clearly,
$E^m \varsubsetneqq    \widetilde{E}^{m,p}$
,
$S^m \varsubsetneqq    \widetilde{S}^{m,p}$
and
$\widetilde{E}^{m,p} \varsubsetneqq    \widetilde{S}^{m,p}$.

It's easy to check that
the
$d_p$ distance, $p\geq 1$, (see \eqref{dpm})
is also a metric on
$\widetilde{S}^{m,p}$.

\section{The characterization of precompact sets in $(\widetilde{S}^{m,p}, d_p)$ }

\begin{df}
   \cite{ma}
   A set $U\subset \widetilde{S}^{m,p}$
   is said to be uniformly $p$-mean bounded if there is a constant $M>0$
   such that
   $d_p(u, 0)\leq M$
   for all $u\in U$.

   We can see that $U$ is uniformly $p$-mean bounded is equivalent to $U$ is a
   bounded set in $( \widetilde{S}^{m,p}, d_p)$.

\end{df}

\begin{df}
  \cite{da}
  Let $u\in \widetilde{S}^{m,p}$. If for given $\varepsilon>0$, there is a $\delta(u,\varepsilon)>0$
  such that for all $0\leq h <\delta$
  $$ \left( \int_h^1  H([u]_\alpha,   [u]_{\alpha-h})^p \;d\alpha   \right)^{1/p}  <\varepsilon ,   $$
where  $1\leq p <+\infty$, then we say $u$ is $p$-mean left-continuous.

Suppose that $U$ is a nonempty set in $\widetilde{S}^{m,p}$.
If the above inequality holds uniformly for all $u\in U$,
then
we say $U$ is  $p$-mean equi-left-continuous.
\end{df}

\begin{tm}
 $U$ is a precompact set in $(\widetilde{S}^{m,p},   d_p)$
 if and only if
 \\
 (\romannumeral1) \  $U$ is uniformly $p$-mean bounded.
 \\
(\romannumeral2) \ $U$ is $p$-mean equi-left-continuous.
\end{tm}

\end{document}